\documentclass{article}

\usepackage{eepic}
\input{epsf.sty}

\DeclareSymbolFont{AMSb}{U}{msb}{m}{n}
\DeclareMathSymbol{\N}{\mathbin}{AMSb}{"4E}
\DeclareMathSymbol{\Z}{\mathbin}{AMSb}{"5A}
\DeclareMathSymbol{\R}{\mathbin}{AMSb}{"52}
\DeclareMathSymbol{\Q}{\mathbin}{AMSb}{"51}
\DeclareMathSymbol{\I}{\mathbin}{AMSb}{"49}
\DeclareMathSymbol{\C}{\mathbin}{AMSb}{"43}

\newcommand{\co}{\supseteq}
\newcommand{\union}{\cup}

\newtheorem{definition}{Definition}

\newtheorem{proposition}[definition]{Proposition}
\newtheorem{corollary}[definition]{Corollary}

\newcommand{\qed}{\hfill \rule{1ex}{1ex}} 

\newenvironment{proof}{{\bf Proof}: }{\qed}

\title{The fine structure of $321$ avoiding permutations.}
\author{M.~H.~Albert}

\begin{document}

\maketitle

\begin{abstract} 
Bivariate generating functions for various subsets of the class of permutations
containing no descending sequence of length three or more are
determined. The notion of absolute indecomposability of a permutation
is introduced, and used in enumerating permutations which have a block
structure avoiding $321$, and whose blocks also have such structure
(recursively). Generalizations of these results are discussed.
\end{abstract}

\section{Introduction}

\begin{figure}
\label{132Diagram}
\centering
\setlength{\unitlength}{0.0005in}
\begingroup\makeatletter\ifx\SetFigFont\undefined%
\gdef\SetFigFont#1#2#3#4#5{%
  \reset@font\fontsize{#1}{#2pt}%
  \fontfamily{#3}\fontseries{#4}\fontshape{#5}%
  \selectfont}%
\fi\endgroup%
{\newcommand{\dashlinestretch}{30}
\begin{picture}(4824,2739)(0,-10)
\path(2862,1062)(4512,1062)(4512,312)
	(2862,312)(2862,1062)
\put(2950,612){\makebox(0,0)[lb]{\smash{{{\SetFigFont{12}{14.4}{\sfdefault}{\mddefault}{\updefault}Avoids 132}}}}}
\path(312,1962)(1962,1962)(1962,1212)
	(312,1212)(312,1962)
\put(400,1512){\makebox(0,0)[lb]{\smash{{{\SetFigFont{12}{14.4}{\sfdefault}{\mddefault}{\updefault}Avoids 132}}}}}
\put(2412,2412){\blacken\ellipse{300}{300}}
\put(2412,2412){\ellipse{300}{300}}
\path(12,2712)(4812,2712)(4812,12)
	(12,12)(12,2712)
\end{picture}
}
\caption{Schematic representation of a $132$ avoiding permutation.}
\end{figure}
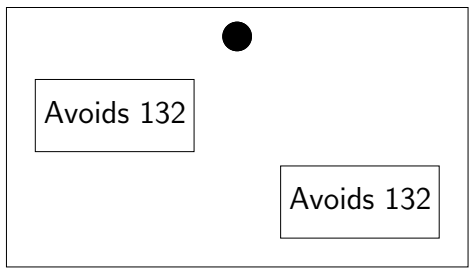

The association of a permutation of $\{1,2,\ldots,n \}$ with its graph
provides a geometric viewpoint in which to consider pattern
avoidance. Thus, for example the permutation avoids $132$ if the
points to the left of the highest point, all lie above the points to
the right, and this condition is true recursively of the points to the
left and to the right of the highest point. This condition is
illustrated in Figure \ref{132Diagram}.  Likewise any set of $n$
points in the plane, no two of which lie on a horizontal or vertical
line can be associated with a permutation of $\{1,2,\ldots,n \}$. In
that case we can represent the geometric information about $132$
avoidance as a picture, with a point for the maximum, and two bounding
rectangles to its left and right, the former lying above the latter,
together with an implicit understanding that the structure within the
rectangles is to be similar.  These associations provide a geometrical
context in which to consider pattern avoidance, and are a common tool
in understanding the class of permutations that avoid one or more
patterns. We explore some ramifications of this viewpoint in the very
simple situation of $321$ avoiding permutations (which we denote by
$A(321)$).

The main purpose of this paper then, is to show how the geometric
context provides a simple method to obtain more detailed enumeration
results about $321$ avoiding permutations than have hitherto been
available. Moreover, these results are obtained uniformly in some
sense. The underlying technique consists of identifying a suitable
geometric configuration which must be attained or avoided, and then
using the structural constraints which that implies in order to
compute the generating function, often multivariate, of the associated
collection of permutations. This technique has applications beyond the
scope of $321$ avoidance, a few of which we consider in the final
sections of the paper.

\begin{figure}
\label{BasicDiagram}
\centering
\setlength{\unitlength}{0.0005in}
\begingroup\makeatletter\ifx\SetFigFont\undefined%
\gdef\SetFigFont#1#2#3#4#5{%
  \reset@font\fontsize{#1}{#2pt}%
  \fontfamily{#3}\fontseries{#4}\fontshape{#5}%
  \selectfont}%
\fi\endgroup%
{\newcommand{\dashlinestretch}{30}
\begin{picture}(2520,2534)(0,-10)
\put(108,612){\blacken\ellipse{200}{200}}
\put(108,612){\ellipse{200}{200}}
\put(708,1212){\blacken\ellipse{200}{200}}
\put(708,1212){\ellipse{200}{200}}
\put(1308,1812){\blacken\ellipse{200}{200}}
\put(1308,1812){\ellipse{200}{200}}
\put(1908,2412){\blacken\ellipse{200}{200}}
\put(1908,2412){\ellipse{200}{200}}
\thicklines
\put(2358,2262){\ellipse{100}{100}}
\put(783,162){\ellipse{100}{100}}
\put(1008,312){\ellipse{100}{100}}
\put(1758,1662){\ellipse{100}{100}}
\put(1608,462){\ellipse{100}{100}}
\put(2208,2112){\ellipse{100}{100}}
\thinlines
\path(1908,2412)(2508,2412)
\path(1308,1812)(2508,1812)
\path(708,1212)(2508,1212)
\path(108,612)(2508,612)
\path(108,612)(108,12)
\path(708,1212)(708,12)
\path(1308,1812)(1308,12)
\path(1908,2412)(1908,12)
\path(2508,2412)(2508,12)
\path(108,12)(2508,12)
\end{picture}
}
\caption{A $321$-avoiding permutation with four left to right
maxima. The given occupancies of the cells define
the permutation ${\bf 4 \, 5} \, 1 \, 2 \, {\bf 7} \, 3 \, 6 \, {\bf
10} \, 8 \, 9 $. Note that the relative horizontal position of
the occupants of the two vertically aligned cells is determined by
the fact that the resulting permutation is not to contain $321$.}
\end{figure}
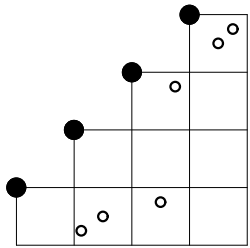

Thought of as a set of points, a permutation avoids $321$ if it does
not contain three points, every pair of which determines a line
segment of negative slope. Of course it is also the case that any
$321$ avoiding permutation is the merge of two increasing
sequences. It is easy to see that one of these sequences can be taken
as the sequence of left to right maxima, that is, those elements which
dominate all of their predecessors.  These elements now determine a
sequence of vertical and horizontal ranges in which the remaining
elements of the permutation must lie. The situation is illustrated in
Figure \ref{BasicDiagram}. Subject to having a fixed number of left to
right maxima, the possible $321$ avoiding permutations that can be
formed are in one to one correspondence with assignments of
non-negative integers to the cells of that diagram, so that no two
cells whose centres are connected by a segment of negative slope are
assigned positive labels. We refer to this diagram, with the cells
containing positive labels simply marked (but the labels themselves
suppressed) as the {\em skeleton} of the permutation.

As is well known, the total number of 132 avoiding permutations of
length $n$ and the total number of 321 avoiding permutations of length
$n$ are the same, both being equal to the $n$th Catalan
number. However, the schematic representation of the $132$ avoiding
permutations makes the correspondence between them and plane binary
trees clear and hence also the equation satisfied by the generating
function of the class, while the corresponding diagram for the $321$
avoiding permutations does not. In some sense $A(321)$ is a class
which exhibits more subtle structure than $A(132)$ does.

Alternatively we might replace each cell in Figure \ref{BasicDiagram}
with its central point, producing a triangular subset of the integer
grid. Given a $321$ avoiding permutation, the points corresponding to
cells with positive labels determine a path in this grid, lying on or
below a diagonal line and
having steps which are of the form $(a,b)$ where $a$ and $b$ are
non-negative integers, not both $0$. This path is also shown in Figure
\ref{BasicDiagram}. Such paths are enumerated in
\cite{Su01} and \cite{Su02}, the former of which also provides a
bijective interpretation of the relationship between their enumeration
and that of standard Delannoy and Schr\"{o}der paths.

Our objective is to find the enumeration via generating functions of
various subsets of the $321$-avoiding permutations. Generally we will
aim to obtain a multivariate generating function (in fact, at most
bivariate) say in $x$ and $y$ where the coefficient of $x^n y^k$ might
be the number of $321$-avoiding permutations of a certain type having
$n$ left to right maxima and $k$ other elements. Then simple
substitutions will allow us to compute either the associated
univariate generating functions, or permutations of size $n$ having
$k$ left to right maxima, or other similar variations. Whenever
possible, the name of a generating function will reflect the name of
the class that it enumerates, so that for example the univariate
generating function for the class $A(321)$ will be $A(t)$, while a
bivariate form might be $A(x,y)$.

To carry out the enumerations of various structurally defined subsets
of the collection of all $321$ avoiding permutations we will
make use of the structural relationships that hold between the subsets
and the class as a whole, and the corresponding
algebraic relationships which hold for the generating
functions. That is, our methodology is firmly in the school
represented by \cite{GJ01}, or \cite{FS01}.

In the next section we rederive equations
describing the basic generating functions for $321$ avoiding
permutations by defining a context-free language whose elements are in
one to one correspondence with $321$ avoiding permutations. This
method is used because it is no longer than extending the results of
\cite{Su01} to the bivariate case, and also because it represents a
technique of wider applicability in the field of pattern class enumeration.

The following section then applies the results to the problem of
enumerating the subsets of $321$ avoiding permutations consisting of:
plus irreducible, minus irreducible, plus indecomposable, and
absolutely irreducible permutations. These terms, and their
significance for enumeration questions are defined below. In the
univariate case, only the last of these is definitely new, although we
have not found detailed expositions of the others in the literature.

In the penultimate section we make use of the enumeration of the absolutely
irreducible $321$ avoiding permutations to enumerate another class,
built from $321$ avoiding permutations by a recursive construction
based upon the wreath product introduced in \cite{AS01}.

In the final section we try to foreshadow future applications of the
methods illustrated here, and mention some connections with other work
in the area of pattern classes.

\section{Enumeration of marked paths}

In Figure \ref{LatticePathDiagram} we rotate the triangular grid from
Figure \ref{BasicDiagram} clockwise by $45^\circ$, and then reflect it
through the $x$-axis. The paths we then obtain begin from $(0,0)$, and
travel through lattice points (the sum of whose coordinates is even)
along segments of slope lying between $-1$ and $1$ inclusive. We will
introduce a grammar which describes such paths. The purpose of
transforming the original diagram is purely psychological. The grammar
which we describe is the language accepted by a certain deterministic
pushdown automaton, and the lattice path then models the number of
elements held in the stack of this automaton as the word is
processed. Given our overlying motivation of considering pattern
classes this also provides a nice link to the generation of
$312$-avoiding permutations using a stack.

Allowing path segments of arbitrary rational slope in $[-1,1]$ would
obviously be a problem in a finite language. We avoid this problem by
adding marks to the end of each path segment, and then replacing each
segment by a horizontal segment followed by one at an inclination of
$\pm 45^{\circ}$. This is illustrated in Figure
\ref{LatticePathDiagram}.

\newcommand{\ch}{{\sf h}}
\newcommand{\cu}{{\sf u}}
\newcommand{\cd}{{\sf d}}
\newcommand{\cm}{{\sf m}}

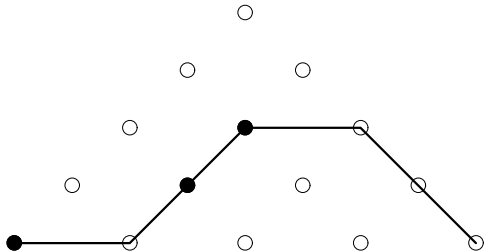
\begin{figure}
\label{LatticePathDiagram}
\centering
\setlength{\unitlength}{0.0005in}
\begingroup\makeatletter\ifx\SetFigFont\undefined%
\gdef\SetFigFont#1#2#3#4#5{%
  \reset@font\fontsize{#1}{#2pt}%
  \fontfamily{#3}\fontseries{#4}\fontshape{#5}%
  \selectfont}%
\fi\endgroup%
{\newcommand{\dashlinestretch}{30}
\begin{picture}(4966,2581)(0,-10)
\put(1283,83){\ellipse{150}{150}}
\put(3683,83){\ellipse{150}{150}}
\put(2483,83){\ellipse{150}{150}}
\put(4883,83){\ellipse{150}{150}}
\put(683,683){\ellipse{150}{150}}
\put(1283,1283){\ellipse{150}{150}}
\put(1883,1883){\ellipse{150}{150}}
\put(2483,2483){\ellipse{150}{150}}
\put(3083,1883){\ellipse{150}{150}}
\put(3683,1283){\ellipse{150}{150}}
\put(4283,683){\ellipse{150}{150}}
\put(3083,683){\ellipse{150}{150}}
\put(83,83){\blacken\ellipse{150}{150}}
\put(83,83){\ellipse{150}{150}}
\put(2483,1283){\blacken\ellipse{150}{150}}
\put(2483,1283){\ellipse{150}{150}}
\put(1883,683){\blacken\ellipse{150}{150}}
\put(1883,683){\ellipse{150}{150}}
\thicklines
\path(83,83)(1283,83)(2483,1283)
	(3683,1283)(4883,83)
\end{picture}
}
\caption{An illustration of the correspondence between occupied
cells in  the triangular grid and lattice paths. Points in the subset
are marked as solid circles, and the corresponding path is
illustrated. The associated word is \cm\ch\cu\cm\cu\cm\ch\cd\cd.}
\end{figure}

We now describe a context free grammar which generates a language
that describes all, and only, the lattice paths which correspond to
marked subsets of the triangular grid of size $n$ for some
$n$. Constructing this grammar is relatively straightforward. We use
four terminal symbols, three of which stand for unit segments in the
path, and one of which represents a vertex in the set:
\ch \ for {\bf h}orizontal segments, 
\cu \ for {\bf u}pward segments, 
\cd \ for {\bf d}ownward segments, and
\cm \ for {\bf m}arking a vertex of the original path.

Informally words in the language are described as follows:
\begin{itemize}
\item
Any change of direction after a \cd \ or a \cu \ requires an intervening
\cm.
\item
\ch's may only occur in blocks immediately following an \cm \ or at the
beginning of a word.
\item
In any initial segment there must be at least as many \cu's as \cd's,
but in the whole word the total number of each is the same.
\item
\cm's can occur anywhere.
\end{itemize}

This description is easily formatted into a grammar, modelled on the
standard grammar corresponding to Dyck paths (which do not allow
horizontal steps, and do not require marks). Each non-terminal
symbol in the grammar represents an excursion, that is a path starting and
ending at the same level and not passing below that level. These
excursions are sometimes restricted by the immediately preceding symbol.
\begin{eqnarray*}
S &\rightarrow& \epsilon \mid \ch S \mid \cu U \cd D \mid \cm M \\
U &\rightarrow& \cm M \mid \cu U \cd D \\
D &\rightarrow& \epsilon \mid \cm M \\
M &\rightarrow& \epsilon \mid \ch S \mid \cu U \cd D
\end{eqnarray*}

As in \cite{CS01} each non-terminal symbol of the grammar is
associated with a generating function (denoted by the same symbol) in
variables $h, u, d, m$. This generating function is obtained by taking
the sum of the monomials corresponding to words represented by that
non-terminal. 

The grammar above is clearly unambiguous since in each rule the
initial symbols of differing productions differ from one another.  So,
it is a simple matter to obtain a system of equations satisfied by the
generating functions of the non-terminals, namely:
\begin{eqnarray*}
S &=& 1 + h S + u d U D + m M \\
U &=& m M + u d U D \\
D &=& 1 + m M \\
M &=& 1 + h S + u d U D.
\end{eqnarray*}

Using a symbolic algebra package or, in a pinch, by hand, this system
can be solved. We are interested principally in the function
$S$ describing words of the language, and this is described as $S =
(1+m) S_1$ where:
\begin{eqnarray*}
\left (-ud{m}^{2}-udm+udmh+ud{m}^{2}h
\right )S_1^{2} &+&  \\
\left (1-h-mh-ud+udh+udmh\right) S_1 -1+u
d && = 0.
\end{eqnarray*}

Given a word in the language, the number of left-to-right maxima in
the $321$ avoiding permutation of which it is the skeleton is equal to
one more than the sum of the number of \ch's and the number of
\cu's. Also, the size of the set that it encodes is equal to the
number of \cm's. So, we can reduce to the generating function $S(x,y)$
where the coefficient of $x^n y^k$ is the number of skeletons with $n$
left to right maxima, and $k$ internal marked cells through the
following substitutions:
\[
h \rightarrow x, \:
u \rightarrow x, \:
d \rightarrow 1, \:
m \rightarrow y,
\]
followed by multiplication by $x$ (and addition of 1 for the empty
graph). This yields:
\begin{eqnarray*}
S(x, y) &=& 1 + x (1 + y) S_1 (x, y) \\
0 &=& (x y + x y^2) S_1^2 + (x y + x - 1) S_1 + 1.
\end{eqnarray*}
The first of these equation can be solved for $S_1$, with the result
being substituted in the second. After some further simplification
this yields:
\begin{equation}
\label{eqS}
x S^2 + (x y + x - 2 y - 1) S + 1 + y = 0.
\end{equation}
By substituting $y = 1$ we will obtain the total number of allowed
markings. So, defining
$S(x) = S(x,1)$:
\[
S(x)^2  + (2x - 3) S(x) + 2 = 0. 
\]
The discriminant of the latter equation is $\sqrt{4 x^2 - 12 x + 1}$
illustrating a connection between these numbers and the Schr\"{o}der
numbers (sequence A001003 of \cite{SP01}). In fact:
\[
S(x) = 1 + 2 x + \sum_{n=1}^{\infty} 2^{n+1} s_{n} x^n
\]
where $s_n$ is the $n$th Schr\"{o}der number. 

This sequence of coefficients also arises as the number of
{\em non-crossing graphs}, that is, graphs with $n$ vertices arranged
as the vertices of a convex polygon, with straight edges connecting
these vertices subject to the condition that no two edges should
intersect at an interior point. This result is due to \cite{DB01},
with a more modern derivation, as well as other related results, given in
\cite{FN01}. 

If we consider the ${n+1 \choose 2}$ cells of the original grid as the
vertices of a graph, $T_n$, two vertices being adjacent if they are
connected by a line of negative slope, then the coefficient of $x^n$
in $S(x)$ counts the independent subsets of $T_n$.  The number of
non-crossing graphs is also the number of independent sets in a
graph. Namely, take as vertices of the graph the possible edges two
such vertices being adjacent if the segments which they represent meet
internally. In this graph $NC_n$, a non-crossing graph corresponds to
an independent set.

So, the generating function for independent subsets of the sequence of
graphs $T_n$ and $NC_{n+1}$ are the same. In fact, inspection of the
results in \cite{FN01} together with a little algebra shows that this
is also true of the bivariate generating functions which mark the
sizes of the independent subsets. That is:

\begin{proposition}
For every $n$ and every $k$, $T_n$ and $NC_{n+1}$ have exactly the
same number of independent subsets of size $k$.
\end{proposition}

However, it is easy to see that for $n \geq 4$, $T_n$ and $NC_{n+1}$
are not isomorphic. For, $T_n$ has exactly four isolated vertices,
while $NC_{n+1}$ has $n+1$ isolated vertices. Detailed expressions for the
coefficients in $S(x)$ and $S(x,y)$ can be found in \cite{FN01}
(Theorem 2, part (ii)) as well as discussions of their asymptotic
expansions.

\section{Consequences for $321$ avoiding permutations}

Before turning to the enumeration of various subsets of the
$321$ avoiding permutations we begin with some remarks about the full
class (whose enumeration is, of course, already well understood
beginning apparently from \cite{Ha01}).

In our original setting, the marked cells arose by considering a
$321$ avoiding permutation having $n$ left-to-right maxima. We argued
that any such permutation corresponded to a labelling of marked
cells with positive integers representing the
number of elements of a permutation contained in a particular
cell. Let $A(x,y)$ be the generating function for $321$ avoiding
permutations where the exponent of $x$ denotes the number of left to
right maxima, and that of $y$ the number of remaining elements. Since
we obtain a $321$ avoiding permutation from its skeleton by replacing
a single cell, marked by a $y$ in $S(x,y)$, by a positive integer, marked
therefore by $y^n$ for some $n > 0$, we obtain:
\[
A(x, y) = S(x, y + y^2 + \cdots) = S \left(x, \, \frac{y}{1-y} \right)
\]
We can also make this substitution in the equation that $S$ satisfies
and then simplify to obtain:
\begin{equation}
\label{eqA1}
y A^2 + (x - y - 1) A + 1 = 0.
\end{equation}
On the other hand, it is perhaps more natural to count permutations of a
common size. So, using $A_{am}$ to denote the generating
function where the coefficient of $x^n y^k$ is the number of $321$
avoiding permutations of length $n$ having $k$ left to right maxima, we obtain:
\[
A_{am}(x, y) = S \left(xy, \, \frac{x}{1-x} \right).
\]
By algebraic manipulation this
function also satisfies a quadratic equation with coefficients
polynomial in $x$ and $y$ namely:
\begin{equation}
\label{eqA2}
x A_{am}^2 + (xy - x - 1) A_{am} + 1 = 0
\end{equation}
A further reduction in complexity occurs when we substitute $y = 1$ in
$A_{am}$ (or $y = x$ in (\ref{eqA1})) giving:
\[
x A(x)^2 - A(x) + 1 = 0
\]
thus confirming, in a rather roundabout way, that the total number of
$321$ avoiding permutations of length $n$ is enumerated by the Catalan
numbers.

The coefficient of $x^n y^k$ in  $A_{am}(x,y)$, which is non-zero only for
$1 \leq k \leq n$ is a {\em Narayana number}, 
\[
[x^n y^k] A = \frac{1}{n} {n \choose k} {n \choose k-1}.
\]
These numbers also arise in \cite{Su01}, but not as a direct
translation of this result since we are no longer in the context of
path counting. They also arise in a number of other contexts including
the enumeration of $k$-way trees (\cite{AW01}) and as the number of
non-crossing partitions of $n$ (\cite{FN01}).

We now turn to the enumeration of various subsets of $A(321)$. First
let us define those classes and the symbols used to specify them:

\begin{definition}
Let $\pi$ be a permutation (in $A(321)$). Then:
\begin{description}
\item[($A_{+irr}$)]
$\pi$ is {\em plus irreducible} if
it does not contain a subword of the form $i \, (i+1)$,
\item[($A_{-irr}$)]
$\pi$ is {\em minus irreducible} if it does not contain a
subword of the form $i \, (i-1)$,
\item[($A_{+ind}$)]
$\pi$ is {\em plus indecomposable} if it does not have a proper
initial segment whose values form an initial segment of $[1,n]$,
\item[($A_{-ind}$)]
$\pi$ is {\em minus indecomposable} if it does not have a proper
final segment whose values form an initial segment of $[1,n]$,
\item[($A_{irr}$)] $\pi$ is {\em absolutely irreducible} if it does
not have a proper subword of length greater than 1 whose values form
an interval in $[1,n]$.
\end{description}
\end{definition}

The irreducible or indecomposable elements of a collection of
permutations can (under suitable closure properties) be thought of as
components in the construction of the other elements of that
class. Again, granted certain closure and uniqueness assumptions, this
can allow enumeration of the entire set based on an enumeration of one
of the collections of components, or vice versa. This particular
exposition of a general combinatorial theme is explored in
\cite{AS01}. We note that the results in that paper could be used to
derive the univariate generating function for the plus irreducibles
and plus indecomposables in $A(321)$ (results which we will rederive
here as a result of obtaining the bivariate form). Furthermore, the only
minus decomposable permutations that avoid $321$ are of the form:
\[
(k+1) \, (k+2) \cdots n \, 1 \, 2 \cdots k
\]
so we will not concern ourselves with that case.

The condition of absolute irreducibility is a new one, and we will
see its application in the next section. The definition is not so
unnatural as it might appear to be at first sight. In terms of the
graph of a permutation it says that if some proper, non-singleton,
part of the permutation is bounded by a rectangle, then there must be
at least one element of the permutation outside of the rectangle but
in either the vertical strip or the horizontal strip determined by
it.

Enumeration results in this section generally take equation
(\ref{eqS}) as their starting point. Recall that this provides the
generating function $S(x,y)$ for skeletons of $321$ avoiding
permutations, with the exponent of $x$ marking the number of left to
right maximals, and that of $y$ the number of occupied cells. So, all
the generating functions we compute will be in the form where the
coefficient of $x^n y^k$ marks the number of permutations of that type
having $n$ left to right maxima and $k$ other elements. As usual, a
simple change of variable, replacing $x$ by $xy$ and $y$ by $x$ would
produce the function enumerating by total number of elements, and
number of left to right maxima.

If $\pi$ is a plus irreducible member of $A(321)$ then no cell can be
occupied by more than one element. Among the diagrams that meet this
criteria, the plus reducible elements contain sequences of more than
one left to right maximum such that the vertical and horizontal bands
which they determine are otherwise empty. Suppose then that we knew
the generating function $A_{+irr}(x,y)$ for the plus irreducible
members of the class. The preceding sentences imply that we would
obtain the generating function $S(x,y)$ by replacing $x$ in
$A_{+irr}(x,y)$ by $x/(1-x)$. So, since the inverse of sending $x$ to
$x/(1-x)$ is to send it to $x/(1+x)$:
\[
A_{+irr}(x,y) = S(x/(1+x), y).
\]
Substitution and simplification in equation \ref{eqS} then yields:
\begin{equation}
\label{eqAPIrr}
x(y+1)A_{+irr}^2 - (xy + 2y + 1)A_{+irr} + (x+1)(y+1) = 0.
\end{equation}
The corresponding univariate form is:
\[
x(x+1)A_{+irr}^2 - (x + 1)^2 A_{+irr} + (x+1)^2 = 0.
\]

An element $\pi$ of $A(321)$ can only be minus reducible if some left to
right maximum $k$ is followed immediately by $k-1$. So
\[
\pi = \alpha k \, (k-1) \, \beta,
\]
for some $\alpha, \beta \in A(321)$ (with $\beta$ of course having all
its values increased by $k$). We can make this decomposition unique by
requiring $k$, $k-1$ to be the first pair of elements witnessing minus
reducibility. Then $\alpha$ is minus irreducible, while $\beta$ could
be any $321$-avoiding permutation. Thus we obtain:
\[
A(x,y) = A_{-irr}(x,y) + A_{-irr}(x,y)(xy)A(x,y).
\]
Or, solving for $A_{-irr}(x,y)$:
\begin{equation}
\label{eqAMIrr}
A_{-irr}(x,y) = \frac{A(x,y)}{1 + xy A(x,y)}.
\end{equation}
The bivariate algebraic equation for $A_{-irr}$ is not very pretty,
but the univariate form is more presentable:
\[
(x^4 + x^2 + x) A_{-irr}^2 + (1 - 2x^2) A_{-irr} + 1 = 0.
\]

Enumerating plus indecomposables is easier  and standard. Every
element of $A(321)$ is either of length $0$ or of the form $\alpha_1
\alpha_2 \cdots \alpha_c$ where each $\alpha_i$ is a plus
indecomposable, shifted upwards by the sum of the lengths of the
preceding $\alpha$'s. Since this decomposition is unique, then using
$A_{+ind}$ to enumerate the non empty plus indecomposables, we obtain:
\[
A = \frac{1}{1 - A_{+ind}},
\]
which can then be readily solved for $A_{+ind}$.

Finally we come to absolute irreducibility. Since the absolutely
irreducibles form a subset of the collection of plus
indecomposables, and of the plus irreducibles, we begin with the form
of the skeleton function which is like that for plus
indecomposables. This already reduces us to permutations that are plus
indecomposable, and plus irreducible in their non left to right maxima.
\[
S_{+ind} (x, y) = \frac{S(x,y)}{1 + S(x,y)}.
\]
Which, by now standard manipulations, satisfies:
\[
(1+y) S_{+ind}^2 - (1 + x + xy) S_{+ind} + xy + x = 0.
\]

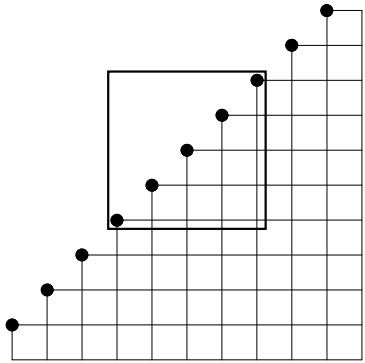
\begin{figure}
\label{BadRectangleDiagram}
\centering
\setlength{\unitlength}{0.0003in}
\begingroup\makeatletter\ifx\SetFigFont\undefined%
\gdef\SetFigFont#1#2#3#4#5{%
  \reset@font\fontsize{#1}{#2pt}%
  \fontfamily{#3}\fontseries{#4}\fontshape{#5}%
  \selectfont}%
\fi\endgroup%
{\newcommand{\dashlinestretch}{30}
\begin{picture}(6120,6134)(0,-10)
\put(108,612){\blacken\ellipse{200}{200}}
\put(108,612){\ellipse{200}{200}}
\put(708,1212){\blacken\ellipse{200}{200}}
\put(708,1212){\ellipse{200}{200}}
\put(1308,1812){\blacken\ellipse{200}{200}}
\put(1308,1812){\ellipse{200}{200}}
\put(1908,2412){\blacken\ellipse{200}{200}}
\put(1908,2412){\ellipse{200}{200}}
\put(2508,3012){\blacken\ellipse{200}{200}}
\put(2508,3012){\ellipse{200}{200}}
\put(3108,3612){\blacken\ellipse{200}{200}}
\put(3108,3612){\ellipse{200}{200}}
\put(3708,4212){\blacken\ellipse{200}{200}}
\put(3708,4212){\ellipse{200}{200}}
\put(4308,4812){\blacken\ellipse{200}{200}}
\put(4308,4812){\ellipse{200}{200}}
\put(4908,5412){\blacken\ellipse{200}{200}}
\put(4908,5412){\ellipse{200}{200}}
\put(5508,6012){\blacken\ellipse{200}{200}}
\put(5508,6012){\ellipse{200}{200}}
\path(1308,1812)(6108,1812)
\path(1908,2412)(6108,2412)
\path(2508,3012)(6108,3012)
\path(3108,3612)(6108,3612)
\path(3708,4212)(6108,4212)
\path(4308,4812)(6108,4812)
\path(4908,5412)(6108,5412)
\path(5508,6012)(6108,6012)
\path(708,1212)(6108,1212)
\path(108,612)(6108,612)
\path(108,12)(6108,12)
\path(4308,4812)(4308,12)
\path(3708,4212)(3708,12)
\path(3108,3612)(3108,12)
\path(2508,3012)(2508,12)
\path(1908,2412)(1908,12)
\path(1308,1812)(1308,12)
\path(708,1212)(708,12)
\path(108,612)(108,12)
\path(4908,5412)(4908,12)
\path(5508,6012)(5508,12)
\path(6108,6012)(6108,12)
\thicklines
\path(1758,2262)(4458,2262)(4458,4962)
	(1758,4962)(1758,2262)
\end{picture}
}
\caption{A potential rectangle for the violation of indecomposability.}
\end{figure}

Consider which non empty rectangles in the diagram associated to an
element of $S_{+ind}$ might not contain other elements inside the
vertical and horizontal strip which they define. In order for this to
hold, the top edge of the rectangle cannot cross a vertical line in
the triangular grid of cells, nor can the left edge cross such a
horizontal line. So, the upper right and lower left corners lie
outside of the grid. Such a rectangle is illustrated in Figure
\ref{BadRectangleDiagram}. The vertical area above the rectangle is
automatically empty as is the horizontal area to the left. So problems
can occur only when we have a non-empty sequence of left to right maxima such
that there are no marked cells in the horizontal or vertical strip
which they define.

If we knew the function $A_{irr}(x,y)$ how could we compute
$S_{+ind}(x,y)$? An element of the latter class could be obtained
beginning from an element of $A_{irr}(x,y)$ by inflating some of the
left to right maxima into a sequence of such maxima, adding no
additional elements in the horizontal or vertical strips which they
determine. If we imagine in Figure \ref{BadRectangleDiagram} that the
illustrated rectangle (and subrectangles of it) are the only ones
which cause a violation of absolute indecomposability, then that
permutation has been constructed by inflating the left to right
maximum just to the left of the rectangle into six such maxima. As
we've already insisted on plus indecomposability, a rectangle whose
leftmost boundary is to the left of the first maximal cannot be
problematic, and so there always is an available maximal to inflate.

There is just a single exception. The permutation $1$ is
absolutely indecomposable, but when we inflate it we do not obtain
plus indecomposable permutations.

Thus, beginning from $A_{irr}(x,y) - x$, we should replace $x$ by
$x/(1-x)$ in order to obtain $S_{+ind}(x,y)$.  Inverting this
replacement we get:
\[
A_{irr} (x,y) = S_{+ind} \left( \frac{x}{1+x}, y \right) - \frac{x}{1+x}.
\]
Carrying out these substitutions and manipulations on the equations
satisfied by the generating function yields:
\begin{equation}
\label{eqAind}
(1+x)(1+y) A_{irr}^2  + (xy - 1) A_{irr} + xy = 0.
\end{equation}
Or in univariate form:
\[
(x+1)^2 A_{irr}^2 + (x^2 - 1) A_{irr} + x^2 = 0.
\]
In fact, this does not quite get us {\em all} the irreducibles as
it omits 1, 12, and the empty permutation. Adjusting the equations to include
this one comes at considerable cost to their appearance, so we prefer
to leave the equation as it stands, adding the necessary $1 + x + x^2$ to the
generating function {\em post facto}. Table \ref{GFCoeffs} summarizes
the sizes of these subsets of $A(321)$.

\begin{table}
\label{GFCoeffs}
\centering
\begin{tabular}{l|rrrrrrrrrrr}
$n$ & 0 & 1 & 2 & 3 & 4 & 5 & 6 & 7 & 8 & 9 & 10 \\ \hline
all & 1 & 1 & 2 & 5 & 14 & 42 & 132 & 429 & 1430 & 4862 & 16796 \\
plus irr. & 1 & 1 & 1 & 2 & 4 & 9 & 21 & 51 & 127 & 323 & 835
\\
minus irr. & 1 & 1 & 1 & 3 & 10 & 31 & 98 & 321 & 1078 & 3686 &
12789 \\
abs. irr. & 1 & 1 & 2 & 0 & 2 & 2 & 7 & 14 & 37 & 90 &
233
\end{tabular}
\caption{Sizes of $A(321)$ and some of its subsets}.
\end{table}

\section{A ``fractal'' class}

As an application of the final results of the preceding section we
will show how the knowledge of the generating function for $A_{irr}$
can be used to compute that of a much more complicated class. At the
risk of further abusing a term which has suffered much abuse already
we would like to introduce the class $F(321)$ of {\em fractal
$321$-avoiders}. These are permutations which, from a distance, appear
to avoid $321$ but which on closer inspection are made up of blocks,
arranged in a $321$-avoiding pattern where each block appears to avoid
$321$ but perhaps on closer inspection is in fact made up out of
blocks \ldots

That is, $\pi \in F(321)$ if either, $\pi \in A(321)$, or $\pi =
\alpha_1 \alpha_2 \cdots \alpha_c$ where
\begin{itemize}
\item
the values occurring in each $\alpha_i$ form an interval,
\item
the permutation represented by $\alpha_i$ is in $F(321)$, and
\item
the relative ordering of the $\alpha_i$, interpreted as a permutation
of length $c$ is in $A(321)$.
\end{itemize}

Many well-known permutation classes can be defined as fractal classes
in this way, or occasionally as natural subclasses of such fractal
classes. For example, the class of separable permutations is precisely
the fractal class generated from the finite base class $\{12, 21\}$.

There is a complementary bottom up description of $F(321)$. Namely,
this class is the closure of the class consisting just of $1$ under
the operation of replacing an element of a permutation by a
$321$-avoiding block. For example:
\[
1 \to 2413 \to 423615
\]
where we initially replace $1$ by $2413$ and then replace the element
$2$ by the permutation $312$ while retaining its relative order within
the entire permutation. Geometrically, we begin with the graph of
$2413$ and then expand the vertex representing $2$ into a copy of the
graph of $312$. Such replacements could just as easily be applied to
each element of a permutation and, in some sense, they already have
been, only $1$ has been replaced by $1$ in three instances. Thus the
two descriptions are equivalent -- the permutation $423615$ consists
of blocks $(423)(6)(1)(5)$ whose relative order is $321$-avoiding, and
where each block is in $F(321)$ (in this instance, in fact in
$A(321)$).

We could also define $F(321)$ algebraically using the wreath product
operator of \cite{AS01} as the smallest non-empty class $X$ satisfying
the equation $X = X \wr A$ where $A = A(321)$. This corresponds to the
bottom up description, while the top down one would suggest $X = A \wr
X$. Consider the first equational description of $X$. Since $X$
contains $A$ we also get that $X$ contains $A \wr A$. Then also $X$
contains $(A \wr A) \wr A$ and so on. Letting $A^n = A^{n-1} \wr A$
for $n > 1$, and $A^1 = A$ we obtain
\[
X \co \union_{n=0}^{\infty} A^n.
\]
On the other hand, the right hand side is contained in its wreath
product with $A$, and so by the definition of $X$:
\[
X = \union_{n=0}^{\infty} A^n.
\]
The second equational definition can be manipulated in the same way
and in fact, as the wreath product is associative, leads to the same
equation, thus confirming that the two approaches are indeed equivalent.

Such an algebraic representation suggests that we ought to be able to
transfer our knowledge of generating functions for $A$ to similar
knowledge about $F$. There is though, a small complication. This
arises from the fact that the choice of blocks to witness the fact
that a permutation belongs to $F(321)$ is not uniquely defined. We
need to obtain uniqueness of some sort if we hope to carry out the
enumeration, and the following general result helps to provide that.

\begin{definition}
Let $\theta$ be a permutation of length $k$ and let $\pi$ be a
permutation of length $n$. Then $\pi$ is {\em
$\theta$-decomposable} if $\pi = \alpha_1 \alpha_2 \cdots \alpha_k$
for some non-empty subwords $\alpha_i$ such that the set of values
occurring in each of the $\alpha_i$ forms an interval and the relative
ordering of these values agrees with the relative ordering of the
corresponding elements of $\theta$. The factorization  $\pi = \alpha_1
\alpha_2 \cdots \alpha_k$ is called a {\em $\theta$-decomposition} of $\pi$.
\end{definition}

With this new definition, we see that a permutation of length 3 or
more is plus decomposable if and only if it is $12$-decomposable,
while a permutation $\pi$ is absolutely irreducible if and only if
it is not $\theta$-decomposable for any $\theta \neq \pi$.

\begin{proposition}
\label{IndecomposablePartition}
Let $\pi$ be an arbitrary permutation. Then there is a unique
absolutely irreducible permutation $\theta$ such that $\pi$ is
$\theta$-decomposable. Moreover, if $\theta \neq 12$ and $\theta
\neq 21$ then the $\theta$-decomposition of $\pi$ is also unique.
\end{proposition}

For example, for $423615$ this decomposition is $(423)(6)(1)(5)$ with
relative ordering $2413$, while for $724513986$ it is $(7) (24513)
(98) (6)$ with relative ordering $3142$. On the other hand $123$ which
is $12$-decomposable admits two such decompositions.

\begin{proof}
Let $\pi$ be given, say $\pi = p_1 p_2 \cdots p_n$. To each $p_i$
associate a maximal proper subword $\alpha_i$ of $\pi$ such that the
values occurring in $\alpha_i$ form an interval (of course, $\alpha_i$
might well be a singleton). 

Suppose that there are $i < j$ such that $\alpha_i$ and $\alpha_j$ overlap
properly but are not equal. Then the elements of $\pi$ belonging to
either $\alpha_i$ or $\alpha_j$ form a subword $\alpha$ whose values
are the union of two overlapping intervals, hence an interval. By the
maximality of either $\alpha_i$ or $\alpha_j$ it must be the case that
$\alpha = \pi$. Then $\alpha_i$ and $\alpha_j$ with the elements
common to $\alpha_i$ deleted form either a $12$ or a $21$
decomposition of $\pi$. These cases are clearly mutually exclusive.

Henceforth suppose that $\pi$ is neither $12$-decomposable nor
$21$-decomposable. Then the $\alpha_i$ form a partition of $\pi$
(i.e. any two are either equal or disjoint). The relative ordering of
the $\alpha_i$ must be some absolutely irreducible permutation
$\theta$, for otherwise we could pool some proper subset of the
$\alpha$'s to form a coarser partition, contradicting the choice of
each $\alpha_i$ as a maximal proper subword of $\pi$ whose values form
an interval. Now reindex the distinct $\alpha_i$ and write $\pi =
\alpha_1 \alpha_2 \cdots \alpha_k$.

Thus we have established the existence of a decomposition of the type
claimed. To establish uniqueness, suppose that another decomposition
of the same kind, say $\pi = \beta_1 \beta_2 \cdots \beta_m$ were
given. We include here the assumption that the relative ordering of
$\beta_1$ through $\beta_m$ forms an absolutely irreducible
permutation. If $\beta_1 \neq \alpha_1$ then $\beta_1$ is a subword of
$\alpha_1$ by the maximality of $\alpha_1$. Now take the least $j$
such that $\beta_1 \beta_2 \cdots \beta_j$ contains $\alpha_1$. Then
in fact we must have $\beta_1 \beta_2 \cdots \beta_j = \alpha_1$ for
otherwise the values in $\alpha_1$ and $\beta_j$ form overlapping
intervals, and so the values in $\beta_1 \beta_2 \cdots \beta_j$ form
an interval, contradicting the maximality of $\alpha_1$ ($\beta_1
\beta_2 \cdots \beta_j \neq \pi$ since $\pi$ is neither $12$ nor
$21$-decomposable). However, $\beta_1 \beta_2 \cdots \beta_j = \alpha_1$
contradicts the absolute irreducibility of the relative ordering of
the $\beta$'s. So $\alpha_1 = \beta_1$. But now the same argument
implies that $\alpha_2 = \beta_2$ and, inductively that in fact $m =
n$ and $\alpha_i = \beta_i$ for all $i$.
\end{proof}

We now return to the analysis of $F(321)$. Let $\pi \in F(321)$ be
given. Suppose that it is neither $12$-decomposable nor
$21$-decomposable. By the proposition above, $\pi = \alpha_1 \alpha_2
\cdots \alpha_k$ for some subwords $\alpha_i$ whose values form
intervals, and whose relative ordering forms an absolutely
irreducible permutation. Since $\pi$ has some decomposition into
subwords whose relative ordering avoids $321$, and since the proof of
the proposition above shows that the $\alpha$'s form the coarsest
possible proper partition of $\pi$ into subwords whose values form
intervals, it must be the case that the relative ordering of the
$\alpha_i$ avoids $321$. Of course, we also have that each $\alpha_i$
belongs to $F(321)$. Conversely, given $\alpha_i$ in $F(321)$, shifted
to have relative order equal to some absolutely irreducible element
$\theta$ of $A(321)$, then, by the very definition of $F(321)$, the
permutation $\alpha_1 \alpha_2 \cdots \alpha_k$ belongs to $F(321)$.

Let $F(x)$ be the generating function for $F(321)$, taken to have
constant term $0$, and $A_{i}(t)$ be the univariate generating
function for the absolutely indecomposable members of $A(321)$ of
length greater than or equal to 3. From the above we see that
$A_{i}(F(x))$ is the generating function for the elements of $F$ of
length at least 3 which are neither $12$-decomposable nor
$21$-decomposable.  Let $F_{+}$ denote the generating function for the
$12$-indecomposable elements of $F$, and $F_{-}$ that of the
$21$-indecomposable elements of $F$, again taken with constant term
$0$. Then $F_{+}F$ enumerates the $12$-decomposable members of $F$
while $F_{-} F$ enumerates the $21$-decomposables. Further relations
arise from the observation that a $12$-indecomposable is either 
$21$-decomposable or both $12$- and $21$-indecomposable, and similarly for
minus indecomposables. We thereby obtain the system of equations:
\begin{eqnarray*}
F &=& x + F_{+}F + F_{-}F + A_{irr}(F) \\
F_{+} &=& x + F_{-}F + A_{irr}(F) \\
F_{-} &=& x + F_{+}F + A_{irr}(F).
\end{eqnarray*}
Solving this system for $F$ gives:
\begin{equation}
\label{eqFrac}
F^2 + (A_{irr}(F) - 1 + x) F + A_{irr}(F) + x = 0.
\end{equation}

We can use the work of the previous section to obtain a radical
expansion of $A_{irr}$:
\[
A_{ind}(x) = {\frac {1-x-\sqrt {-3\,{x}^{2}-2\,x+1}}{2(x+
1)}} - x^2.
\]
Then substitution in (\ref{eqFrac}) and elimination of radicals gives:
\begin{eqnarray}
{F}^{6}+\left (-2\,x+3\right ){F}^{4}+\left (-2\,x-1\right ){F}^{3}+
 && \nonumber \\ 
\label{eqFrac2}
\left (-3\,x+3+{x}^{2}\right ){F}^{2} + \left (2\,{x}^{2}-1-2\,x\right )
F+x+{x}^{2} &=& 0.
\end{eqnarray}
The first few terms of the associated power series are:
\[
\begin{array}{l}
x+2\,{x}^{2}+6\,{x}^{3}+24\,{x}^{4}+116\,{x}^{5}+625\,{x}^{6}+3580\,{
x}^{7}+  \\
21297\,{x}^{8}+130084\,{x}^{9}+810737\,{x}^{10}+O\left ({x}^{11
}\right )
\end{array}
\]
and the exponential constant governing the growth rate, is the
reciprocal of the radius of convergence of this series. This radius is
the least positive root of the discriminant of (\ref{eqFrac2}), which is
an irreducible polynomial of degree $7$.
The value of the exponential constant is approximately $7.346751$,
compared to $4$ for the underlying class $A(321)$.

\section{Summary and Conclusions}

We began this paper with an explicitly constructed grammar to describe
the skeletons of $321$ avoiding permutations. In general there is a
close connection between combinatorial classes with algebraic
generating functions and unambiguous context free languages. This
connection can either be used, as here, to provide an explicit
enumeration of a class, or to provide a ``soft'' proof that the
generating function of a class is in fact algebraic. The former
approach has become much more attractive with the ready availability
of symbolic algebra packages since the algebraic manipulations
necessary to solve the equations arising from the grammar are
undeniably tedious. The latter approach has been used in \cite{AA01}
to provide algebraicity results for a family of pattern classes. It
can also be used in the context of generating functions for
generating trees, thereby generalizing a number of the theorems in
\cite{BBDFGG01} about the existence of algebraic generating
functions. It must be noted though that the results of that paper and
similar results in \cite{BF01} provide much more explicit detail
concerning the generating functions that they produce.

One of the striking features of the equations for the various
irreducible and indecomposable subsets of $A(321)$ is their
simplicity. In some sense then the enumerative coincidences that we
observed are not so startling, since there is a relatively limited
supply of simple quadratic equations. Because of the simple form of
these equations, one could apply Lagrange inversion to obtain explicit
formulae for many of the coefficients, as is done for example in
\cite{FN01}, or indeed carry out detailed asymptotic analyses of these
coefficients. 

We restricted ourselves to bivariate generating functions but the
reader should note that the techniques employed can be naturally
applied to produce other statistics of these permutations. For
example, it would be a simple matter to produce, if one wished, a
generating function $A(x,y,z,w)$ where $x$ marked total size, $y$
marked left to right maxima, $z$ marked the number of occurences of $i
(i+1)$ among the left to right maxima, and $w$ that number among the
remaining elements. 

The class $A(321)$ is the simplest pattern class, in terms of the
patterns which it avoids, that contains infinitely many absolutely
irreducibles. The techniques used in the preceding section to solve
(in the sense of enumeration) the wreath fixed point equation:
\[
X = A \wr X
\]
apply, owing to proposition \ref{IndecomposablePartition}, completely
generally to any base class $A$ in which the absolutely
irreducibles can be enumerated. In particular for a positive
integer $n$ let $D_n$ be the class of permutations which ``fractally
have $\leq n$ elements''. That is, they are comprised of at most $n$
blocks, each of which is comprised of at most $n$ blocks, each of
which \ldots Then $D_n$ is the solution of the fixed point equation:
\[
X = F_n \wr X
\]
where $F_n$ is the class of permutations of size $\leq n$. Since $F_n$
is finite equation (\ref{eqFrac}) is simply a polynomial and we obtain:

\begin{corollary}
Each of the classes $D_n$ has an algebraic generating function.
\end{corollary}

There is much further information to be gleaned from the
representation of a class as a subclass of $D_n$ when this is
possible, and we hope to explore these matters in a future paper.

One aspect of $F(321)$ that has been notably omitted is a
description in terms of minimal forbidden patterns. It appears that
this set may be finite consisting of:
\[
25314, \, 35142, \, 41352, \, 42513, \,  362514, \, 531642
\]
but all that can be said with certainty at this point is that no
further minimal forbidden patterns exist of length 15 or less.

\section{Acknowledgements}

Many of the connections between these problems and known results would
have been missed without the electronic version of the encyclopedia of
integer sequences (\cite{S01}). Mike Atkinson listened patiently to
many preliminary, and incorrect, expositions, and also provided a key
step in the proof of Proposition \ref{IndecomposablePartition}. {\em
Maple} did most of the hard work.

\bibliographystyle{plain}

\begin{thebibliography}{10}

\bibitem{AA01}
M.H. Albert and M.D. Atkinson.
\newblock Sorting with a forklift.
\newblock {\em Otago University CS Tech. Report}, 2002-06:25pp, 2002.

\bibitem{AW01}
M.~D. Atkinson and Louise Walker.
\newblock Enumerating {$k$}-way trees.
\newblock {\em Inform. Process. Lett.}, 48(2):73--75, 1993.

\bibitem{AS01}
M.D. Atkinson and T.~Stitt.
\newblock Restricted permutations and the wreath product.
\newblock {\em Discrete Math.}, 259:19--36, 2002.

\bibitem{BBDFGG01}
Cyril Banderier, Mireille Bousquet-M{\'e}lou, Alain Denise, Philippe Flajolet,
  Dani{\`e}le Gardy, and Dominique Gouyou-Beauchamps.
\newblock Generating functions for generating trees.
\newblock {\em Discrete Math.}, 246(1-3):29--55, 2002.


\bibitem{BF01}
Cyril Banderier and Philippe Flajolet.
\newblock Basic analytic combinatorics of directed lattice paths.
\newblock {\em Theoret. Comput. Sci.}, 281(1-2):37--80, 2002.


\bibitem{CS01}
N.~Chomsky and M.~P. Sch{\"u}tzenberger.
\newblock The algebraic theory of context-free languages.
\newblock In {\em Computer programming and formal systems}, pages 118--161.
  North-Holland, Amsterdam, 1963.

\bibitem{DB01}
C.~Domb and A.~J. Barrett.
\newblock Enumeration of ladder graphs.
\newblock {\em Discrete Math.}, 9:341--358, 1974.

\bibitem{FS01}
P.~Flajolet and R.~Sedgewick.
\newblock {\em Analytic Combinatorics---Symbolic Combinatorics}.
\newblock Preprint, http://algo.inria.fr/flajolet/Publications/books.html,
  2002.

\bibitem{FN01}
Philippe Flajolet and Marc Noy.
\newblock Analytic combinatorics of non-crossing configurations.
\newblock {\em Discrete Math.}, 204(1-3):203--229, 1999.

\bibitem{GJ01}
I.~P. Goulden and D.~M. Jackson.
\newblock {\em Combinatorial enumeration}.
\newblock A Wiley-Interscience Publication. John Wiley \& Sons Inc., New York,
  1983.

\bibitem{Ha01}
J.~M. Hammersley.
\newblock A few seedlings of research.
\newblock In {\em Proceedings of the Sixth Berkeley Symposium on Mathematical
  Statistics and Probability (Univ. California, Berkeley, Calif., 1970/1971),
  Vol. I: Theory of statistics}, pages 345--394, Berkeley, Calif., 1972. Univ.
  California Press.

\bibitem{S01} N.~J.~A. Sloane, editor (2002), The On-Line Encyclopedia
of Integer Sequences, published electronically at \\
\texttt{http://www.research.att.com/~njas/sequences/}.


\bibitem{SP01}
N.~J.~A. Sloane and Simon Plouffe.
\newblock {\em The encyclopedia of integer sequences}.
\newblock Academic Press Inc., San Diego, CA, 1995.


\bibitem{Su02}
Robert~A. Sulanke.
\newblock A symmetric variation of a distribution of {K}reweras and {P}oupard.
\newblock {\em J. Statist. Plann. Inference}, 34(2):291--303, 1993.

\bibitem{Su01}
Robert~A. Sulanke.
\newblock Counting lattice paths by {N}arayana polynomials.
\newblock {\em Electron. J. Combin.}, 7(1):Research Paper 40, 9 pp.
  (electronic), 2000.

\end{thebibliography}

\end{document}